\documentclass[11pt]{amsart}

\usepackage[T1]{fontenc}
\usepackage{lmodern}
\usepackage{microtype}
\usepackage{amsmath,amssymb,amsthm}
\usepackage{booktabs}
\usepackage{array}
\usepackage[margin=1.15in]{geometry}
\usepackage{xcolor}
\usepackage[colorlinks=true,linkcolor=blue!60!black,citecolor=blue!60!black,urlcolor=blue!60!black]{hyperref}

\numberwithin{equation}{section}

\theoremstyle{plain}
\newtheorem{theorem}{Theorem}[section]
\newtheorem{lemma}[theorem]{Lemma}

\newtheorem{corollary}[theorem]{Corollary}
\theoremstyle{remark}
\newtheorem{remark}[theorem]{Remark}

\newcommand{\s}{\mathsf{s}}
\newcommand{\F}{\mathbb{F}}
\newcommand{\Z}{\mathbb{Z}}
\newcommand{\PP}{\mathbb{P}}
\newcommand{\E}{\mathbb{E}}
\newcommand{\cP}{\mathcal{P}}
\newcommand{\cW}{\mathcal{W}}
\newcommand{\bits}[1]{\texttt{#1}}

\title[Exponential lower bounds for Erd\H{o}s--Ginzburg--Ziv constants]{Cyclic fibers and exponential lower bounds for\\ Erd\H{o}s--Ginzburg--Ziv constants}
\author{Krishnendu Bhowmick}
\address{Krishnendu Bhowmick \\ TU Chemnitz, Germany.}
\email{krish.combi@gmail.com}
\date{}


\begin{document}

\begin{abstract}
We show that for every integer $n\ge2$ that is not a power of two there
are constants $c_n>0$ and $B_n>2$ such that $\s(C_n^r)\ge c_nB_n^r$ for
all $r\ge1$. For powers of two, Harborth's formula gives growth rate exactly
$2$. For the remaining even exponents, the previous lower bounds were of
order $r2^r$. When $n=2m$ with $m\ge3$ odd, we prove the transfer
inequality
\[
A_n(H\oplus C_n^5)\ \ge\ 12\,A_n(H\oplus C_n)+8\,A_n(H\oplus C_m),
\]
where $A_n(G)$ is the maximal length of a sequence over $G$ without an $n$-term
zero-sum subsequence and $H$ is any finite abelian group whose exponent
divides $n$. The transfer extends a configuration of Elsholtz to twenty full
cyclic fibers. Combined with a nine-point binary construction for
$2m$-term zero-sums over $C_m^r$, it gives bases at least $2.0341$ that tend to
$2.0394\ldots$; combined with a finite-field construction, it gives
$\s(C_6^r)\ge c\cdot 2.0491^r$. For exponents divisible by $4$ that are not
powers of two, we obtain a base above two from a weighted random construction
on the eighteen-point configuration of Wang and Zeng.
\end{abstract}

\maketitle

\section{Introduction}

Let $G$ be a finite abelian group and let $L$ be a positive integer divisible
by the exponent $\exp(G)$. We write $\s_L(G)$ for the least integer $s$ such
that every sequence of $s$ elements of $G$ has a subsequence of exactly $L$
terms with sum zero. Repeated elements are allowed, and subsequences are
chosen by position. The \emph{Erd\H{o}s--Ginzburg--Ziv constant} is
$\s(G)=\s_{\exp(G)}(G)$. It is often more convenient to work with
\[
A_L(G)=\s_L(G)-1,
\]
the largest length of a sequence over $G$ that has no zero-sum subsequence of
length $L$. We study how $\s(C_n^r)$ grows with the rank $r$ when the exponent
$n\ge2$ is fixed and $r\ge1$. In intermediate generalized constants we also
allow rank zero, with $C_n^0=\{0\}$ and $A_L(\{0\})=L-1$.

Harborth~\cite{Har73} showed that $\s(C_n^r)\ge 2^r(n-1)+1$ for all $n$, with
equality when $n$ is a power of two. For odd $n\ge3$, Elsholtz~\cite{Els04}
improved this to
\begin{equation}\label{eq:odd}
\s(C_n^r)\ \ge\ (n-1)\,9^{\lfloor r/3\rfloor}2^{\,r-3\lfloor r/3\rfloor}+1,
\end{equation}
whose exponential rate is $9^{1/3}>2$. For even $n$ that is not a power of
two, Wang and Zeng~\cite{WZ26} showed that Harborth's bound is not sharp once
$r\ge4$, and Elsholtz~\cite{Els26} showed that the excess over Harborth's bound
grows at least like $c_n\,r\,2^r$. These lower bounds still have exponential
base two. We obtain a base strictly larger than two for every such exponent.

\begin{theorem}\label{thm:main}
For every integer $n\ge 2$,
\[
\liminf_{r\to\infty}\s(C_n^r)^{1/r}
\begin{cases}
=2, & \text{if $n$ is a power of two},\\
>2, & \text{otherwise}.
\end{cases}
\]
\end{theorem}

The even exponents that are not powers of two split into two cases. When
$n\equiv 2\pmod 4$ we obtain explicit bases from a transfer inequality between
generalized and classical constants (Theorem~\ref{thm:transfer}).

\begin{theorem}\label{thm:twice-odd}
Let $m\ge3$ be odd, and define $\rho_m>1$ and $\alpha_m>2$ by
\[
\rho_m=\left(\frac{9^{2m}}{9+36\binom{2m}{m}}\right)^{1/(5(2m-1))},
\qquad
\alpha_m^5=12\alpha_m+8\rho_m .
\]
Then, for every $r\ge1$,
\begin{equation}\label{eq:twice-odd}
\s_{2m}(C_m^r)\ \ge\ \frac{m}{12e}\,\rho_m^r+1,
\qquad
\s(C_{2m}^r)\ \ge\ \frac{m}{12e}\Bigl(\frac{2}{\alpha_m}\Bigr)^{4}\alpha_m^r+1 .
\end{equation}
Moreover, $\rho_m\ge\rho_3=3^{6/25}=1.3016945937\ldots$ and hence
$\alpha_m\ge\alpha_3=2.0341017978\ldots$ for all odd $m\ge3$, and
\[
\lim_{\substack{m\to\infty\\ m\text{ odd}}}\rho_m=(9/2)^{1/5}=1.3509600385\ldots,
\qquad
\lim_{\substack{m\to\infty\\ m\text{ odd}}}\alpha_m=2.0394244419\ldots .
\]
\end{theorem}

For $n=6$ a finite-field construction gives a better input to the transfer.

\begin{theorem}\label{thm:six}
Let $\alpha_*>2$ be the root of $\alpha_*^5=12\alpha_*+8\sqrt[3]{3}$, so
$\alpha_*=2.0491202034\ldots$. Then, for every $r\ge1$,
\[
\s(C_6^r)\ \ge\ c\,\alpha_*^r+1,\qquad c=2\cdot3^{-5/3}\Bigl(\frac{2}{\alpha_*}\Bigr)^4 .
\]
\end{theorem}

For the remaining even exponents that are not powers of two, we use a
weighted random construction on the configuration of Wang and Zeng. For
$0<p<1$ and $n=Km$ put
\begin{equation}\label{eq:Q}
Q=Q_{n,m}(p)=\sum_{j=0}^{K}\binom{n}{jm}p^{jm}(1-p)^{n-jm},
\qquad
Z=14+2Q^{-1/(n-1)} .
\end{equation}
Thus $Q$ is the probability that a binomial random variable with parameters
$n$ and $p$ is divisible by $m$.

\begin{theorem}\label{thm:general}
Let $n=Km$, where $K=2^a$ with $a\ge1$ and $m\ge3$ is odd, and let $0<p<1$.
Then $Z>16$. If $r=4t+d\ge1$ with $t\ge0$ and $0\le d\le3$, then
\begin{equation}\label{eq:general}
\s(C_n^r)\ \ge\ \frac{n}{4e}\,2^dZ^t+1 .
\end{equation}
In particular, with $B_{n,p}=Z^{1/4}>2$,
\[
\s(C_n^r)\ \ge\ \frac{n}{4e}\Bigl(\frac{2}{B_{n,p}}\Bigr)^3B_{n,p}^r+1
\qquad\text{for every } r\ge1 .
\]
\end{theorem}

Theorem~\ref{thm:general} also covers $n\equiv2\pmod4$, but there it gives
weaker bases than Theorem~\ref{thm:twice-odd}; see
Remark~\ref{rem:twice-odd-weaker}. The improvement is in the exponential
base as $r\to\infty$; for small ranks, the earlier constructions may give
larger lower bounds.

\begin{proof}[Proof of Theorem~\ref{thm:main}]
If $n$ is a power of two, Harborth's formula $\s(C_n^r)=2^r(n-1)+1$ gives
$\s(C_n^r)^{1/r}\to2$. If $n\ge3$ is odd, \eqref{eq:odd} gives
$\liminf_r\s(C_n^r)^{1/r}\ge9^{1/3}>2$. If $n$ is even but not a power of two,
then $n=2^am$ with $a\ge1$ and $m\ge3$ odd, and Theorem~\ref{thm:general}
gives $\liminf_r\s(C_n^r)^{1/r}\ge B_{n,p}>2$.
\end{proof}

\subsection*{Comparison with generalized constants}
The first inequality in \eqref{eq:twice-odd} concerns zero-sums of $2m$ terms
over $C_m^r$, not $m$ terms. For this problem Bitz, Griffith and
He~\cite{BGH20} obtained a family of lower bounds whose exponential bases
tend to $5/4$ as $m\to\infty$. Our bases tend to $(9/2)^{1/5}>5/4$ as $m$
grows through odd integers. This compares the two families in the
large-modulus limit, rather than asserting optimality for each fixed $m$.
Geneson~\cite{Gen17} improved their bound for zero-sums of length $km$ with
$k>2$, a different range of prescribed lengths.

\subsection*{Organization}
Section~\ref{sec:tools} contains two standard tools: the alteration method and
binary padding. Section~\ref{sec:transfer} proves the transfer inequality.
Section~\ref{sec:binary} gives the binary construction for $\s_{2m}(C_m^r)$,
and Section~\ref{sec:field} gives the finite-field construction and proves
Theorems~\ref{thm:twice-odd} and~\ref{thm:six}. These sections are
self-contained. Section~\ref{sec:general} proves Theorem~\ref{thm:general}; it
relies on a classification of zero-sums due to Elsholtz~\cite{Els26}.

\section{Two tools}\label{sec:tools}

The first tool is the alteration method: sample a random sequence, then delete
one term from each forbidden subsequence. For prescribed zero-sum lengths this
was used by Bitz, Griffith and He~\cite{BGH20}, and for higher-degree
zero-sum constants by Costa and Della Fiore~\cite[Theorem~2.5]{CDF24}. We need
it for nonuniform distributions.

\begin{lemma}\label{lem:alteration}
Let $L\ge2$, let $G$ be a finite abelian group with $\exp(G)\mid L$, and let
$\mu$ be a probability distribution on $G$. Suppose that the sum of $L$
independent samples from $\mu$ is zero with probability $\theta>0$. Put
\[
D_L=\bigl((L-1)!\bigr)^{1/(L-1)},\qquad N=\bigl\lfloor D_L\,\theta^{-1/(L-1)}\bigr\rfloor .
\]
Then $A_L(G)\ge(1-1/L)N$. If $L\ge4$, then also
\[
A_L(G)\ \ge\ \frac{L}{4e}\,\theta^{-1/(L-1)} .
\]
\end{lemma}

\begin{proof}
Choose $N$ independent elements with distribution $\mu$, keeping their
positions distinct even when their values coincide. Let $Y$ be the number of
$L$-element sets of positions whose values sum to zero. Since
$N^{L-1}\theta\le D_L^{L-1}=(L-1)!$,
\[
\E Y=\binom{N}{L}\theta\le\frac{N^L\theta}{L!}\le\frac{N}{L}.
\]
Fix a realization with $Y\le N/L$ and delete one position from each of the $Y$
forbidden sets. Deleting terms creates no new forbidden set, so the remaining
sequence has no zero-sum subsequence of length $L$, and its length is at least
$N-N/L$.

For the second assertion, comparing the sum with an integral gives
\[
\log\bigl((L-1)!\bigr)=\sum_{k=2}^{L-1}\log k\ \ge\ \int_1^{L-1}\log x\,dx\ \ge\ (L-1)\log(L-1)-(L-1),
\]
so $D_L\ge(L-1)/e$. Since $D_L\ge1$ and $\theta\le1$, the number
$D_L\theta^{-1/(L-1)}$ is at least $1$, and $\lfloor x\rfloor\ge x/2$ for
$x\ge1$. Hence
\[
A_L(G)\ \ge\ \Bigl(1-\frac1L\Bigr)\frac{L-1}{2e}\,\theta^{-1/(L-1)}
=\frac{(L-1)^2}{2eL}\,\theta^{-1/(L-1)}\ \ge\ \frac{L}{4e}\,\theta^{-1/(L-1)},
\]
where the last step uses $2(L-1)^2\ge L^2$, which holds for $L\ge4$.
\end{proof}

The second tool lifts a sequence to a group with more cyclic coordinates.

\begin{lemma}\label{lem:padding}
Let $n\ge2$, let $G$ and $H$ be finite abelian groups whose exponents divide
$n$, and let $d\ge0$. Then
\[
A_n(G\oplus C_n^d)\ \ge\ 2^d A_n(G)
\qquad\text{and}\qquad
A_n(G\oplus H)\ \ge\ A_n(G).
\]
\end{lemma}

\begin{proof}
Let $S$ be a sequence over $G$ of length $A_n(G)$ with no $n$-term zero-sum.
Replace each term $x$ of $S$ by the $2^d$ terms $(x,b)$, $b\in\{0,1\}^d$. In
an $n$-term zero-sum of the new sequence, each of the last $d$ coordinates is
a sum of $n$ numbers in $\{0,1\}$ that is divisible by $n$, so it equals $0$
or $n$. Hence all selected terms have the same $b$. For that $b$, the selected
terms come from distinct terms of $S$ and form an $n$-term zero-sum in $S$,
which is impossible. The second inequality follows by embedding $G$ as
$G\oplus\{0\}$.
\end{proof}

\section{A transfer through full cyclic fibers}\label{sec:transfer}

Throughout this section $m\ge3$ is odd and $n=2m$.

\begin{theorem}\label{thm:transfer}
For every finite abelian group $H$ whose exponent divides $n$,
\begin{equation}\label{eq:transfer}
A_n(H\oplus C_n^5)\ \ge\ 12\,A_n(H\oplus C_n)+8\,A_n(H\oplus C_m).
\end{equation}
Suppose moreover that $c>0$ and $\beta>1$ satisfy $A_n(C_m^s)\ge c\beta^s$ for
all $s\ge0$. Let $\alpha$ be the unique positive root of
\begin{equation}\label{eq:alpha}
\alpha^5=12\alpha+8\beta .
\end{equation}
Then $\alpha>2$, and for every $r\ge1$,
\begin{equation}\label{eq:transfer-rate}
\s(C_n^r)\ \ge\ c\Bigl(\frac2\alpha\Bigr)^4\alpha^r+1 .
\end{equation}
\end{theorem}

\subsection{The fibers}
By the Chinese remainder theorem we identify $C_n=C_2\oplus C_m$, so that
\[
C_n^5=\F_2^5\oplus C_m^5,\qquad C_m^5=C_m\oplus C_m^4 .
\]
We write elements of $C_m^5$ as $(\tau,w)$ with $\tau\in C_m$ and $w\in C_m^4$.
Let $e_1,\dots,e_4$ be the standard basis of $\Z^4$, and let $V\subset\Z^4$
consist of the twelve \emph{spoke vectors} $\pm e_i\pm e_4$ ($1\le i\le3$) and
the eight \emph{cube vectors} $(\pm1,\pm1,\pm1,0)$. For each $v\in V$,
Table~\ref{tab:labels} gives two labels $\pi^0(v),\pi^1(v)\in\F_2^5$. Put
$\delta_v=\pi^0(v)+\pi^1(v)$ and
\[
F_v=\bigl\{\bigl(\pi^0(v)+(\tau\bmod 2)\,\delta_v,\ (\tau\bmod m,\ v\bmod m)\bigr):\tau\in\Z\bigr\}.
\]
For a spoke vector, $\delta_v\ne0$ and $F_v$ is a coset of the cyclic subgroup
of order $n$ generated by $(\delta_v,(1,0))$. For a cube vector,
$\pi^0(v)=\pi^1(v)$, so $\delta_v=0$ and $F_v$ is a coset of the cyclic
subgroup of order $m$ generated by $(0,(1,0))$. The coordinates of distinct
vectors of $V$ differ by at most $2<m$, so distinct vectors stay distinct
modulo $m$ and the twenty fibers are pairwise disjoint. The points with
$\tau\in\{0,1\}$ are the forty points of~\cite[Table~3]{Els26}; we use the
whole fibers.

\begin{table}[ht]
\centering
\small
\begin{tabular}{>{$}c<{$}cc@{\qquad}>{$}c<{$}cc}
\toprule
v & $\pi^0(v)$ & $\pi^1(v)$ & v & $\pi^0(v)$ & $\pi^1(v)$\\
\midrule
(0,0,1,1)   & \bits{00000} & \bits{01110} & (1,1,1,0)    & \bits{10100} & \bits{10100}\\
(0,0,1,-1)  & \bits{00010} & \bits{01100} & (1,1,-1,0)   & \bits{01101} & \bits{01101}\\
(0,0,-1,1)  & \bits{10111} & \bits{11001} & (1,-1,1,0)   & \bits{00011} & \bits{00011}\\
(0,0,-1,-1) & \bits{11011} & \bits{10101} & (1,-1,-1,0)  & \bits{11010} & \bits{11010}\\
(0,1,0,1)   & \bits{01011} & \bits{00101} & (-1,0,0,1)   & \bits{10001} & \bits{10011}\\
(0,1,0,-1)  & \bits{01001} & \bits{00111} & (-1,0,0,-1)  & \bits{11111} & \bits{11101}\\
(0,-1,0,1)  & \bits{11100} & \bits{11110} & (-1,1,1,0)   & \bits{00001} & \bits{00001}\\
(0,-1,0,-1) & \bits{10000} & \bits{10010} & (-1,1,-1,0)  & \bits{11000} & \bits{11000}\\
(1,0,0,1)   & \bits{01010} & \bits{00110} & (-1,-1,1,0)  & \bits{10110} & \bits{10110}\\
(1,0,0,-1)  & \bits{00100} & \bits{01000} & (-1,-1,-1,0) & \bits{01111} & \bits{01111}\\
\bottomrule
\end{tabular}
\medskip
\caption{The labels $\pi^0(v),\pi^1(v)$, taken from~\cite[Table~3]{Els26}.
Each string lists the five binary coordinates $\pi_1\pi_2\pi_3\pi_4\pi_5$ in
order.}
\label{tab:labels}
\end{table}

\subsection{Separation}

\begin{lemma}\label{lem:separation}
Every $n$-term zero-sum sequence supported on $\bigcup_{v\in V}F_v$ is
supported on a single fiber $F_v$.
\end{lemma}

\begin{proof}
Let $T$ be such a sequence. Each term of $T$ lies in exactly one fiber; for
the $j$-th term let $v^{(j)}\in V$ be its vector, regarded as an integer
vector, and let $\pi^{(j)}\in\{\pi^0(v^{(j)}),\pi^1(v^{(j)})\}$ be its binary
part. Since $T$ has zero sum,
\[
\sum_{j=1}^n v^{(j)}=m\,y\quad\text{for some } y\in\{-2,-1,0,1,2\}^4,
\qquad
\sum_{j=1}^n\pi^{(j)}_k\equiv0\pmod2\quad(1\le k\le5),
\]
where the bound on $y$ holds because each coordinate of the left-hand sum has
absolute value at most $n=2m$.

\emph{Step 1: parity.} Inspection of Table~\ref{tab:labels} shows that, for
every $v\in V$ and for both labels $\pi\in\{\pi^0(v),\pi^1(v)\}$,
\begin{equation}\label{eq:parity}
\begin{aligned}
v_1&\equiv\pi_2+\pi_3+\pi_5 \pmod2, &\qquad v_2&\equiv\pi_1+\pi_5\pmod2,\\
v_3+1&\equiv\pi_1+\pi_2+\pi_3\pmod2, & v_1+v_2+v_3+3&\equiv2\pi_1\pmod4.
\end{aligned}
\end{equation}
Sum each identity over the $n$ terms of $T$. Every binary coordinate has an
even total, so the first two identities give $my_1\equiv my_2\equiv0$ and the
third gives $my_3+n\equiv0\pmod2$. Since $m$ is odd, $y_1,y_2,y_3$ are even.
In the fourth identity, $2\sum_j\pi^{(j)}_1\equiv0\pmod4$, so
$m(y_1+y_2+y_3)+6m\equiv0\pmod4$, and hence $y_1+y_2+y_3\equiv2\pmod4$. As
$y_1,y_2,y_3\in\{-2,0,2\}$, an odd number of them is nonzero.

\emph{Step 2: three nonzero coordinates.} If $y_1,y_2,y_3$ are all nonzero,
then for each $i\le3$ the integers $v^{(j)}_i\in\{-1,0,1\}$ sum to $\pm n$, so
they all equal the same sign $\varepsilon_i$. Only one vector of $V$ has
first three coordinates $(\varepsilon_1,\varepsilon_2,\varepsilon_3)$, namely
a cube vector, so $T$ lies in a single fiber.

\emph{Step 3: one nonzero coordinate.} Suppose instead that exactly one of
them is nonzero, say $y_i=2\varepsilon$ with $\varepsilon\in\{\pm1\}$. As in
Step~2, every $v^{(j)}$ lies on the face
\[
V_{i,\varepsilon}=\{v\in V:v_i=\varepsilon\}
=\{u_+,u_-\}\cup\{\text{four cube vectors}\},\qquad u_\pm=\varepsilon e_i\pm e_4 .
\]
The six vectors of the face form three \emph{opposite pairs} $\{v,w\}$ with
$v+w=2\varepsilon e_i$: the pair $\{u_+,u_-\}$ and two pairs of cube
vectors. Let $k,l$ be the two indices in $\{1,2,3\}\setminus\{i\}$. Since
$y_k=y_l=0$ and $u_\pm$ vanish in coordinates $k$ and $l$, the cube vectors
with $(v_k,v_l)$ equal to $(1,1)$, $(-1,-1)$, $(1,-1)$, $(-1,1)$ occur with multiplicities
$N_{++},N_{--},N_{+-},N_{-+}$ satisfying
\[
N_{++}+N_{+-}-N_{-+}-N_{--}=0,\qquad N_{++}-N_{+-}+N_{-+}-N_{--}=0 .
\]
Hence $N_{++}=N_{--}=:a_1$ and $N_{+-}=N_{-+}=:a_2$, so opposite cube vectors
occur equally often. Let $a_\pm$ be the multiplicity of $u_\pm$. Counting
terms and looking at the fourth coordinate gives
\[
2a_1+2a_2+a_++a_-=2m,\qquad a_+-a_-=my_4 .
\]
The first equation shows that $a_++a_-$, and hence $a_+-a_-$, is even, so
$y_4$ is even. If $y_4=\pm2$, then $a_\pm=2m$ and $T$ lies in the fiber
$F_{u_\pm}$. It remains to rule out $y_4=0$. In that case $a_+=a_-=:a_0$ and
\begin{equation}\label{eq:abc}
a_0+a_1+a_2=m .
\end{equation}

\emph{Step 4: a binary character.} For each face, let $\lambda=\lambda_{i,\varepsilon}$ be the linear
functional on $\F_2^5$ whose coefficient vector is given below:
\[
\begin{array}{c|cccccc}
(i,\varepsilon) & (1,-1) & (1,1) & (2,-1) & (2,1) & (3,-1) & (3,1)\\
\hline
\lambda & \bits{00100} & \bits{00010} & \bits{00100} & \bits{00110} & \bits{00110} & \bits{00110}
\end{array}
\]
Using Table~\ref{tab:labels} one checks two properties on each face:
$\lambda(\delta_v)=0$ for all six $v\in V_{i,\varepsilon}$, and
$\lambda(\pi^0(v))+\lambda(\pi^0(w))=1$ for each of the three opposite pairs
$\{v,w\}$. The first property says that $\lambda$ takes the same value
$\lambda(\pi^0(v))$ on both labels of $v$, and hence on all of $F_v$. Since
the two vectors of an opposite pair occur equally often, the second property
and~\eqref{eq:abc} give
\[
\lambda\Bigl(\sum_{j=1}^n\pi^{(j)}\Bigr)\equiv a_0+a_1+a_2=m\equiv1\pmod2 .
\]
This contradicts $\sum_j\pi^{(j)}=0$.
\end{proof}

\subsection{Proof of Theorem~\ref{thm:transfer}}
For a spoke vector $v$, let $x_v=(\pi^0(v),(0,v\bmod m))$ and
$g_v=(\delta_v,(1,0))$, and define the injective map
\[
\varphi_v\colon H\oplus C_n\to H\oplus C_n^5,\qquad (h,\tau)\mapsto(h,\ x_v+\tau g_v).
\]
For a cube vector $v$, define $\varphi_v\colon H\oplus C_m\to H\oplus C_n^5$
by the same formula with $g_v=(0,(1,0))$. The image of $\varphi_v$ is
$H\oplus F_v$, and these twenty sets are pairwise disjoint.

For each spoke vector choose a sequence of length $A_n(H\oplus C_n)$ over
$H\oplus C_n$ with no $n$-term zero-sum; for each cube vector choose such a
sequence of length $A_n(H\oplus C_m)$ over $H\oplus C_m$. Let $S$ be the
concatenation of their images under the maps $\varphi_v$. Its length is the
right-hand side of~\eqref{eq:transfer}.

Suppose that $S$ has an $n$-term zero-sum subsequence. Its projection to
$C_n^5$ is an $n$-term zero-sum supported on $\bigcup_vF_v$, so by
Lemma~\ref{lem:separation} all its terms lie in a single $H\oplus F_v$. Write
them as $\varphi_v(h_j,\tau_j)$. Since $\exp(H\oplus C_n^5)$ divides $n$, we
have $nx_v=0$, so
\[
\sum_{j=1}^n\varphi_v(h_j,\tau_j)=\Bigl(\sum_j h_j,\ \Bigl(\sum_j\tau_j\Bigr)g_v\Bigr).
\]
As $\tau\mapsto\tau g_v$ is injective, $\sum_jh_j=0$ and $\sum_j\tau_j=0$,
so the terms $(h_j,\tau_j)$ form an $n$-term zero-sum in the sequence chosen
for $v$. This contradiction proves~\eqref{eq:transfer}.

For the rate, put $\Phi(s,u)=A_n(C_n^s\oplus C_m^u)$. Applying
\eqref{eq:transfer} with $H=C_n^{s-5}\oplus C_m^u$ gives, for $s\ge5$,
\[
\Phi(s,u)\ \ge\ 12\,\Phi(s-4,u)+8\,\Phi(s-5,u+1).
\]
We show by induction on $s$, simultaneously for all $u\ge0$, that
\begin{equation}\label{eq:induction}
\Phi(s,u)\ \ge\ c\Bigl(\frac2\alpha\Bigr)^4\alpha^s\beta^u .
\end{equation}
For $0\le s\le4$, Lemma~\ref{lem:padding} gives
$\Phi(s,u)\ge2^sA_n(C_m^u)\ge2^sc\beta^u$, and
$2^s\ge(2/\alpha)^4\alpha^s$ because $\alpha>2$ and $s\le4$. For $s\ge5$, the
recursion and the induction hypothesis give~\eqref{eq:induction}, because
$12\alpha^{-4}+8\beta\alpha^{-5}=1$ by~\eqref{eq:alpha}. Taking $u=0$ and
$s=r$ proves~\eqref{eq:transfer-rate}.

Finally, $x^5-12x-8\beta$ has one sign change in its coefficients, so by
Descartes' rule of signs it has exactly one positive root. It is negative at
$x=2$ because $\beta>1$, so this root exceeds $2$.
\qed

\section{A binary construction for prescribed zero-sum length}\label{sec:binary}

Consider the nine binary vectors
\begin{equation}\label{eq:nine}
\begin{aligned}
P_1&=(0,0,0,0,0), & P_2&=(0,0,1,1,0), & P_3&=(0,0,1,1,1),\\
P_4&=(0,1,0,1,0), & P_5&=(0,1,1,0,0), & P_6&=(0,1,1,1,1),\\
P_7&=(1,0,0,1,1), & P_8&=(1,0,1,0,0), & P_9&=(1,0,1,1,0),
\end{aligned}
\end{equation}
and let $\cP=\{P_1,\dots,P_9\}$. The rows are grouped by their first two
coordinates, which are $00$, $01$ and $10$ respectively.

\begin{lemma}\label{lem:nine}
Let $m\ge2$, and consider a sequence of $2m$ elements of $\cP$ with sum zero
in $C_m^5$. Let $\mu_1,\dots,\mu_9$ be the multiplicities of
$P_1,\dots,P_9$ divided by $m$. Then either every $\mu_i$ is an integer, or
$(\mu_i)$ equals $1/2$ on one of the sets
\[
\{P_1,P_2,P_4,P_5\},\qquad\{P_1,P_3,P_7,P_8\},\qquad\{P_4,P_6,P_7,P_9\}
\]
and vanishes elsewhere. Consequently:
\begin{enumerate}
\item if $m$ is odd, a sequence of $2m$ elements of $\cP$ has sum zero in
$C_m^5$ if and only if it consists of $2m$ copies of one element or of $m$
copies each of two elements;
\item for every $m\ge2$, the same holds for sequences of $2m$ elements of
$\cP'=\cP\setminus\{P_1,P_7\}$.
\end{enumerate}
\end{lemma}

\begin{proof}
The $\mu_i$ are nonnegative and sum to $2$. Since the $P_i$ have entries $0$
and $1$, each coordinate of the integer sum of the sequence is divisible by
$m$, so $\sum_i\mu_iP_i\in\Z^5$. The first two coordinates show that the
group masses
\[
M_0=\mu_1+\mu_2+\mu_3,\qquad M_1=\mu_4+\mu_5+\mu_6,\qquad M_2=\mu_7+\mu_8+\mu_9
\]
are integers, since $M_2$ and $M_1$ are these coordinates and $M_0=2-M_1-M_2$.
Thus either one group has mass $2$, or two groups have mass $1$ each.

\emph{One group.} If only $M_0$ is nonzero, the last three coordinates are
$\mu_2+\mu_3,\ \mu_2+\mu_3,\ \mu_3$; these are integers, so $\mu_3$, $\mu_2$
and $\mu_1=2-\mu_2-\mu_3$ are integers. For $M_1$ the last three coordinates
are $\mu_5+\mu_6,\ \mu_4+\mu_6,\ \mu_6$, and for $M_2$ they are
$\mu_8+\mu_9,\ \mu_7+\mu_9,\ \mu_7$; the same argument applies.

\emph{$M_0=M_1=1$.} Using $\mu_2+\mu_3=1-\mu_1$ and
$\mu_5+\mu_6=1-\mu_4$, the last three coordinates show that
\[
\mu_1+\mu_4,\qquad\mu_1+\mu_5,\qquad\mu_3+\mu_6
\]
are integers. All $\mu_i$ lie in $[0,1]$. If $\mu_1\in\{0,1\}$, then
$\mu_4,\mu_5$ are integers, hence so are $\mu_6=1-\mu_4-\mu_5$, $\mu_3$ and
$\mu_2$. Otherwise $0<\mu_1<1$ forces $\mu_4=\mu_5=1-\mu_1$, so
$\mu_6=2\mu_1-1$ and $\mu_1\ge1/2$. As $0\le\mu_3\le1-\mu_1$,
\[
2\mu_1-1\ \le\ \mu_3+\mu_6\ \le\ \mu_1 .
\]
For $1/2<\mu_1<1$ this interval lies inside $(0,1)$ and contains no integer.
For $\mu_1=1/2$ it forces $\mu_3=\mu_6=0$, so
$\mu_1=\mu_2=\mu_4=\mu_5=1/2$.

\emph{$M_0=M_2=1$.} Now
\[
\mu_1+\mu_7,\qquad\mu_1+\mu_8,\qquad\mu_3+\mu_7
\]
are integers. Again $\mu_1\in\{0,1\}$ forces all $\mu_i$ to be integers.
Otherwise $\mu_7=\mu_8=1-\mu_1$, so $\mu_9=2\mu_1-1$ and $\mu_1\ge1/2$. As
$0\le\mu_3\le1-\mu_1$,
\[
1-\mu_1\ \le\ \mu_3+\mu_7\ \le\ 2-2\mu_1 .
\]
For $1/2<\mu_1<1$ this interval lies inside $(0,1)$ and contains no integer.
For $\mu_1=1/2$ it forces $\mu_3+\mu_7=1$, so $\mu_3=1/2$ and
$\mu_2=\mu_9=0$; thus $\mu_1=\mu_3=\mu_7=\mu_8=1/2$.

\emph{$M_1=M_2=1$.} Now
\[
\mu_4+\mu_7,\qquad\mu_5+\mu_8,\qquad\mu_6+\mu_7
\]
are integers, so $\mu_4-\mu_6\in\{-1,0,1\}$. If $\mu_4-\mu_6=\pm1$, then
$\mu_4,\mu_6\in\{0,1\}$ and all $\mu_i$ are integers. Otherwise
$\mu_4=\mu_6=x$ and $\mu_5=1-2x$ with $0\le x\le1/2$. If $x=0$, all $\mu_i$
are integers. If $0<x<1/2$, then $\mu_7=1-x$ and $\mu_8=2x$, so
$\mu_9=1-\mu_7-\mu_8=-x<0$, which is impossible. If $x=1/2$, then $\mu_5=0$,
$\mu_7=1/2$, $\mu_8=0$ and $\mu_9=1/2$; thus
$\mu_4=\mu_6=\mu_7=\mu_9=1/2$.

This proves the first assertion. For (1), the half-integral cases require
$m/2$ to be an integer, which fails for odd $m$. For (2), each exceptional set
contains $P_1$ or $P_7$. In the remaining integral case the $\mu_i$ are
nonnegative integers summing to $2$, which gives the two stated forms. The
converse holds because $2mP=0$ and $m(P+P')=0$ in $C_m^5$.
\end{proof}

\begin{theorem}\label{thm:binary}
For every odd integer $m\ge3$, put
\begin{equation}\label{eq:qrho}
q_m=\frac{9+36\binom{2m}{m}}{9^{2m}},\qquad\rho_m=q_m^{-1/(5(2m-1))} .
\end{equation}
If $r=5t+u$ with $t\ge1$ and $0\le u\le4$, then
\begin{equation}\label{eq:binary-floor}
A_{2m}(C_m^r)\ \ge\ \Bigl(1-\frac1{2m}\Bigr)\Bigl\lfloor\bigl((2m-1)!\bigr)^{1/(2m-1)}\rho_m^{5t}\Bigr\rfloor .
\end{equation}
Moreover, for every $r\ge0$,
\begin{equation}\label{eq:binary-clean}
A_{2m}(C_m^r)\ \ge\ \frac{m}{12e}\,\rho_m^r\ \ge\ \frac{m}{12e}\,3^{6r/25} .
\end{equation}
The witnessing sequences may be chosen in $\{0,1\}^r$. Finally,
\[
3^{6/25}=\rho_3\ \le\ \rho_m\ \le\ 9^{1/5},\qquad
\lim_{\substack{m\to\infty\\ m\text{ odd}}}\rho_m=(9/2)^{1/5}.
\]
\end{theorem}

\begin{proof}
Let $\mu$ be the uniform distribution on $\cP\subset C_m^5$. By
Lemma~\ref{lem:nine}(1), among the $9^{2m}$ ordered $2m$-tuples from $\cP$
exactly $9$ are constant zero-sums and $\binom92\binom{2m}{m}$ are nonconstant
zero-sums. So the zero-sum probability is $q_m$. On $C_m^{5t}$ we use the
product distribution $\mu^{\otimes t}$; a sum vanishes exactly when each block
of five coordinates does, so the probability is $q_m^t$. Lemma~\ref{lem:alteration}
with $L=2m$ and $\theta=q_m^t$ gives~\eqref{eq:binary-floor} for $r=5t$, and
Lemma~\ref{lem:padding} (adding $u$ zero coordinates) extends it to $r=5t+u$.
All terms lie in $\{0,1\}^r$.

Since the constant tuples alone give $q_m\ge9^{1-2m}$, we have
$\rho_m\le9^{1/5}$, and so $\rho_m^u\le9^{4/5}<6$ for $u\le4$. For $t\ge1$,
the second part of Lemma~\ref{lem:alteration} and padding give
\[
A_{2m}(C_m^r)\ \ge\ \frac{2m}{4e}\,\rho_m^{5t}=\frac{m}{2e}\,\rho_m^{r}\rho_m^{-u}\ \ge\ \frac{m}{12e}\,\rho_m^r .
\]
For $r\le4$, a sequence of $2m-1$ zeros shows that $A_{2m}(C_m^r)\ge2m-1$,
which exceeds $(m/(12e))\cdot6>(m/(12e))\rho_m^r$. This proves the first
inequality in~\eqref{eq:binary-clean}.

For the lower bound on $\rho_m$, note that $9+36\binom{6}{3}=729=3^6$. When
$m$ increases by one, $9+36\binom{2m}{m}$ grows by a factor less than $4$,
since $\binom{2m+2}{m+1}/\binom{2m}{m}=2(2m+1)/(m+1)<4$, whereas
$3^{(8m+6)/5}$ grows by the factor $3^{8/5}>4$. By induction,
$9+36\binom{2m}{m}\le3^{(8m+6)/5}$ for every integer $m\ge3$. Hence
$q_m\le3^{-6(2m-1)/5}$ and $\rho_m\ge3^{6/25}=\rho_3$, which also gives the
second inequality in~\eqref{eq:binary-clean}. Finally,
$\binom{2m}{m}^{1/(2m)}\to2$ gives the limit.
\end{proof}

\begin{remark}\label{rem:nonmonotone}
The bound $\rho_m\ge\rho_3$ does not rely on monotonicity. In fact, the
sequence is not monotone: $\rho_3<(9/2)^{1/5}$, whereas
$\rho_{21}=1.35110\ldots>(9/2)^{1/5}$, and its limit is $(9/2)^{1/5}$.
\end{remark}

Part~(2) of Lemma~\ref{lem:nine} gives a weaker bound with no parity
restriction.

\begin{corollary}\label{cor:seven}
For every integer $m\ge2$, put
\[
\rho'_m=\left(\frac{7^{2m}}{7+21\binom{2m}{m}}\right)^{1/(5(2m-1))}.
\]
Then $A_{2m}(C_m^r)\ge(m/(12e))\,{\rho'_m}^{r}$ for every $r\ge0$, with
witnessing sequences in $\{0,1\}^r$. Moreover, $\rho'_m\to(7/2)^{1/5}$ as
$m\to\infty$.
\end{corollary}

\begin{proof}
Repeat the proof of Theorem~\ref{thm:binary} with the uniform distribution on
the seven points of $\cP'$. By Lemma~\ref{lem:nine}(2) the zero-sum
probability is $(7+21\binom{2m}{m})/7^{2m}$, and $\rho'_m\le7^{1/5}$ gives
${\rho'_m}^u\le7^{4/5}<6$.
\end{proof}

\section{A finite-field input for six-term sums}\label{sec:field}

For exponent six, a finite-field construction gives a stronger input to
Theorem~\ref{thm:transfer} than Theorem~\ref{thm:binary}. Such power
parametrizations are classical in coding theory: Dumer~\cite[Section~III,
Lemma~2 and Theorem~1]{Dum95} uses columns defined by finite-field powers and
controls their small-support linear dependencies through the rank of their
locators. We use the following construction as an auxiliary input and give
a direct proof using Newton's identities and the discriminant.

\begin{lemma}\label{lem:field}
Let $k\ge1$ be odd, let $F=\F_{3^k}$, and put
\[
T=\{(x,x^2,x^4):x\in F\}\subset F^3 .
\]
Every six-term zero-sum sequence supported on $T$ consists either of one
element taken six times or of two distinct elements taken three times each.
Consequently, since $F^3\cong C_3^{3k}$ as groups,
\[
A_6(C_3^{3k})\ \ge\ 2\cdot3^k .
\]
\end{lemma}

\begin{proof}
Suppose that $\sum_{i=1}^6(x_i,x_i^2,x_i^4)=0$. Let $p_j=\sum_ix_i^j$ be the
power sums and $\sigma_j$ the elementary symmetric functions of
$x_1,\dots,x_6$. By assumption $p_1=p_2=p_4=0$. Newton's identities
\[
p_1=\sigma_1,\qquad p_2=\sigma_1p_1-2\sigma_2,\qquad
p_4=\sigma_1p_3-\sigma_2p_2+\sigma_3p_1-4\sigma_4
\]
then give $\sigma_1=\sigma_2=\sigma_4=0$, because $2$ and $4$ are invertible
in characteristic three. Hence
\[
f(X)=\prod_{i=1}^6(X-x_i)=X^6-\sigma_3X^3-\sigma_5X+\sigma_6,
\qquad f'(X)=-\sigma_5 .
\]
If $\sigma_5\ne0$, then $f$ and $f'$ have no common root, so the $x_i$ are
distinct, and
\[
\prod_{i<j}(x_i-x_j)^2=(-1)^{15}\prod_{i=1}^6f'(x_i)=-\sigma_5^6 .
\]
The left-hand side is a nonzero square in $F$. The right-hand side is not a
square, because $3^k\equiv3\pmod4$ for odd $k$, so $-1$ is a nonsquare in
$F$. Hence $\sigma_5=0$.

The Frobenius map $y\mapsto y^3$ is an automorphism of $F$, so there are
$b,c\in F$ with $b^3=-\sigma_3$ and $c^3=\sigma_6$. Then
$f(X)=(X^2+bX+c)^3$. Every root of $f$ therefore has multiplicity divisible by
three, which gives the two stated possibilities for the multiset
$\{x_1,\dots,x_6\}$. The map $x\mapsto(x,x^2,x^4)$ is injective, so the same
holds for the terms of the sequence. In particular, the sequence containing
each point of $T$ twice has no six-term zero-sum, and its length is
$2\cdot3^k$.
\end{proof}

\begin{proof}[Proof of Theorem~\ref{thm:twice-odd}]
By Theorem~\ref{thm:binary}, $\s_{2m}(C_m^r)-1=A_{2m}(C_m^r)\ge(m/(12e))\rho_m^r$
for all $r\ge0$, which is the first bound. Theorem~\ref{thm:transfer} with
$c=m/(12e)$ and $\beta=\rho_m$ gives the second. The root of
$x^5-12x-8\beta$ above two increases with $\beta$, so $\rho_m\ge\rho_3$ gives
$\alpha_m\ge\alpha_3$. The derivative $5x^4-12$ is positive for $x\ge2$, so by
the implicit function theorem this root depends continuously on $\beta$, and
$\rho_m\to(9/2)^{1/5}$ gives the limit of $\alpha_m$.
\end{proof}

\begin{proof}[Proof of Theorem~\ref{thm:six}]
Put $\beta=3^{1/3}$ and $c_0=2\beta^{-5}$. We first show that
$A_6(C_3^s)\ge c_0\beta^s$ for all $s\ge0$. For $s\ge3$, let $k$ be the
largest odd integer with $3k\le s$; then $0\le s-3k\le5$. By
Lemmas~\ref{lem:field} and~\ref{lem:padding},
$A_6(C_3^s)\ge A_6(C_3^{3k})\ge2\cdot3^k=2\beta^{3k}\ge2\beta^{s-5}$. For
$s\in\{0,1,2\}$, a sequence of five zeros gives $A_6(C_3^s)\ge5\ge2\beta^{s-5}$.
Now Theorem~\ref{thm:transfer} with $m=3$ gives
$\s(C_6^r)\ge c_0(2/\alpha_*)^4\alpha_*^r+1$, which is the claim.
\end{proof}

\section{All even exponents}\label{sec:general}

In this section $n=Km$, where $K=2^a$ with $a\ge1$ and $m\ge3$ is odd. We use
the configuration of Wang and Zeng~\cite{WZ26} as presented by
Elsholtz~\cite[Table~2]{Els26}. For each row $R$ of
Table~\ref{tab:wz} and each listed value of $t$, let
\[
g(R,t)=(y,\ x-z,\ q,\ t)\in C_n^4,
\]
where $(x,y,z)$ is the triple in row $R$ and $q\in C_n$ is determined by
$q\equiv x\pmod m$ and $q\equiv b_R(t)\pmod K$. Let $\cW$ be the set of these
eighteen points.

\begin{table}[ht]
\centering
\small
\begin{tabular}{c>{$}c<{$}>{$}c<{$}>{$}c<{$}}
\toprule
row $R$ & (x,y,z) & b_R(t) & \text{values of }t\\
\midrule
A & (0,0,0) & -t & 0,1\\
B & (1,0,0) & -1 & 2,\ 2+K\\
C & (0,1,0) & -1 & 0,1\\
D & (0,0,1) & -1 & 1,\ 1+K\\
E & (1,0,1) & 1-t & 0,1\\
F & (0,1,1) & -1 & 0,1\\
G & (1,1,2) & 0 & 0,1\\
H & (1,2,2) & 0 & 0,1\\
I & (2,1,2) & 0 & 0,1\\
\bottomrule
\end{tabular}
\medskip
\caption{The Wang--Zeng configuration, following~\cite[Table~2]{Els26}.}
\label{tab:wz}
\end{table}

Put $h_1=g(B,2)$, $\ell_1=g(B,2+K)$, $h_2=g(D,1)$ and $\ell_2=g(D,1+K)$. In
rows B and D the value $b_R(t)$ does not depend on $t$, so the third
coordinates of $h_i$ and $\ell_i$ agree and
\begin{equation}\label{eq:pair-diff}
\ell_i-h_i=Ke_4\qquad(i=1,2),
\end{equation}
an element of order $m$.

\begin{lemma}[Elsholtz~{\cite[Proposition~5.1]{Els26}}]\label{lem:wz}
Every $n$-term zero-sum sequence supported on $\cW$ is either constant or
supported on one of the pairs $\{h_1,\ell_1\}$, $\{h_2,\ell_2\}$. A sequence
of $n$ terms supported on $\{h_i,\ell_i\}$ has sum zero exactly when the
number of terms equal to $\ell_i$ is divisible by $m$.
\end{lemma}

The first assertion is~\cite[Proposition~5.1]{Els26}. In~\cite{Els26} the
points of $\cW$ carry prescribed multiplicities ($n-1$ for the sixteen points
$g(R,t)$ with $t\in\{0,1,2\}$ and $m-1$ for $\ell_1,\ell_2$), but the
classification there holds for every sequence supported on $\cW$, with
arbitrary multiplicities, and this is what we need. The second assertion
follows from~\eqref{eq:pair-diff}: if $j$ terms equal $\ell_i$, the sum is
$nh_i+jKe_4=jKe_4$, which vanishes exactly when $m\mid j$.

\begin{proof}[Proof of Theorem~\ref{thm:general}]
Write $\gamma=Q^{-1/(n-1)}$, so that $Z=14+2\gamma$. Define a probability
distribution $\mu$ on $\cW$ as follows. Each of the fourteen points other than
$h_1,\ell_1,h_2,\ell_2$ has probability $1/Z$; for $i=1,2$, the point $h_i$
has probability $\gamma(1-p)/Z$ and $\ell_i$ has probability $\gamma p/Z$.
These probabilities sum to $(14+2\gamma)/Z=1$.

Let $X_1,\dots,X_n$ be independent with distribution $\mu$. By
Lemma~\ref{lem:wz}, the sum vanishes in exactly one of sixteen disjoint ways:
the $X_j$ are all equal to one of the fourteen unpaired points, or they all lie
in $\{h_i,\ell_i\}$ for some $i$ with the number equal to $\ell_i$ divisible by
$m$. (The constant sequences at $h_i$ and $\ell_i$ are included in the second
kind.) The first kind contributes $14Z^{-n}$. For each $i$, the second kind
contributes
\[
\sum_{j=0}^{K}\binom{n}{jm}\Bigl(\frac{\gamma p}{Z}\Bigr)^{jm}\Bigl(\frac{\gamma(1-p)}{Z}\Bigr)^{n-jm}
=\Bigl(\frac{\gamma}{Z}\Bigr)^nQ .
\]
Since $\gamma^{n-1}Q=1$,
\begin{equation}\label{eq:theta}
\PP(X_1+\dots+X_n=0)=\frac{14+2\gamma^nQ}{Z^n}=\frac{14+2\gamma}{Z^n}=Z^{1-n}.
\end{equation}

On $C_n^{4t}$ use the product distribution $\mu^{\otimes t}$. A sum vanishes
exactly when each block of four coordinates does, so by independence
and~\eqref{eq:theta} the zero-sum probability is $Z^{-(n-1)t}$. Since $n\ge6$,
Lemma~\ref{lem:alteration} gives $A_n(C_n^{4t})\ge(n/(4e))Z^t$, and
Lemma~\ref{lem:padding} gives
$A_n(C_n^{4t+d})\ge2^d(n/(4e))Z^t$. This is~\eqref{eq:general}.

Next, $Q>0$ because of the term $j=0$, and $Q<1$ because the binomial
variable equals $1$ with positive probability and $m\nmid1$. Hence
$\gamma>1$, $Z>16$ and $B_{n,p}>2$. Finally, for $r=4t+d$,
\[
2^dZ^t=2^dB_{n,p}^{r-d}=\Bigl(\frac{2}{B_{n,p}}\Bigr)^dB_{n,p}^r\ \ge\ \Bigl(\frac{2}{B_{n,p}}\Bigr)^3B_{n,p}^r .
\qedhere
\]
\end{proof}

\begin{remark}\label{rem:twice-odd-weaker}
For $n=2m$, the fiber construction gives a larger base than
Theorem~\ref{thm:general} for every choice of $p$. Indeed,
\[
Q_{2m,m}(p)=(1-p)^{2m}+\binom{2m}{m}p^m(1-p)^m+p^{2m}
>\Bigl(\frac23\Bigr)^{2m-1}.
\]
To see this, use symmetry to assume $p\le1/2$. If $p\le1/4$, the first
term is at least $(3/4)^{2m}>(2/3)^{2m-1}$. If $1/4\le p\le1/2$, the
middle term is at least $\binom{2m}{m}(3/16)^m$. At $m=3$ this is
$135/1024>32/243$; increasing $m$ by one multiplies this lower bound by
$3(2m+1)/(8(m+1))\ge21/32>4/9$, so the inequality follows by induction.
Consequently
\[
B_{2m,p}<17^{1/4}=2.030543\ldots<\alpha_3\le\alpha_m.
\]
\end{remark}

\begin{remark}
An analogue of Theorem~\ref{thm:transfer} for $n=Km$ with $K\ge4$ would
extend the cyclic-fiber method to the remaining even exponents. To improve
Theorem~\ref{thm:general}, such a transfer would need a recurrence whose
growth rate exceeds $B_{n,p}$.
\end{remark}
\section*{AI disclaimer}
 This work was carried out in collaboration with AI, namely ChatGPT 6 Astra. 

\bibliographystyle{plain}
\bibliography{references}

@article{BGH20,
  author = {Bitz, Jared and Griffith, Sarah and He, Xiaoyu},
  title = {Exponential lower bounds on the generalized {Erd\H{o}s--Ginzburg--Ziv} constant},
  journal = {Discrete Mathematics},
  volume = {343},
  number = {12},
  year = {2020},
  pages = {112083},
  doi = {10.1016/j.disc.2020.112083},
  url = {https://doi.org/10.1016/j.disc.2020.112083}
}

@article{CDF24,
  author = {Costa, Simone and Della Fiore, Stefano},
  title = {Bounds on the higher degree {Erd\H{o}s--Ginzburg--Ziv} constants over {$\mathbb{F}_q^n$}},
  journal = {Archiv der Mathematik},
  volume = {122},
  year = {2024},
  pages = {17--29},
  doi = {10.1007/s00013-023-01916-4},
  url = {https://doi.org/10.1007/s00013-023-01916-4}
}

@misc{Els26,
  author = {Elsholtz, Christian},
  title = {The {Erd\H{o}s--Ginzburg--Ziv} constant of {$C_n^r$} for even {$n$}: rigid sequences, products, and growth beyond {Harborth}'s bound},
  year = {2026},
  howpublished = {arXiv:2610.35819v2},
  url = {https://arxiv.org/abs/2610.35819v2}
}

@article{Els04,
  author = {Elsholtz, Christian},
  title = {Lower bounds for multidimensional zero sums},
  journal = {Combinatorica},
  volume = {24},
  year = {2004},
  pages = {351--358},
  doi = {10.1007/s00493-004-0022-y},
  url = {https://doi.org/10.1007/s00493-004-0022-y}
}

@article{Har73,
  author = {Harborth, Heiko},
  title = {{Ein Extremalproblem f{\"u}r Gitterpunkte}},
  journal = {Journal f{\"u}r die reine und angewandte Mathematik},
  volume = {262/263},
  year = {1973},
  pages = {356--360}
}

@misc{Gen17,
  author = {Geneson, Jesse},
  title = {Improved lower bound on generalized {Erd\H{o}s--Ginzburg--Ziv} constants},
  year = {2017},
  howpublished = {arXiv:1712.02069},
  url = {https://arxiv.org/abs/1712.02069}
}

@article{Dum95,
  author = {Dumer, Ilya},
  title = {Nonbinary double-error-correcting codes designed by means of algebraic varieties},
  journal = {IEEE Transactions on Information Theory},
  volume = {41},
  number = {6},
  year = {1995},
  pages = {1657--1666},
  doi = {10.1109/18.476238},
  url = {https://doi.org/10.1109/18.476238}
}

@misc{WZ26,
  author = {Wang, Guoqing and Zeng, Xiangneng},
  title = {Inverse zero-sum problems of {Erd\H{o}s--Ginzburg--Ziv} constant: {Property C} and {Property D} conjectures},
  year = {2026},
  howpublished = {arXiv:2610.25271v2},
  url = {https://arxiv.org/abs/2610.25271v2}
}

\end{document}